\documentclass[12pt,reqno,a4paper]{amsart}
\usepackage{amsfonts,amsmath,amssymb}
\usepackage{algorithmic,algorithm}
\usepackage{amsthm,array,a4wide,graphics}

\usepackage{hyperref}
\usepackage{color}
\definecolor{blau}{rgb}{0,0,0.75} %color for in-document links
\hypersetup{colorlinks,linkcolor=blau,citecolor=blue}

\allowdisplaybreaks

\newcommand{\shuffk}{\ensuremath{\scriptscriptstyle{\boldsymbol{\sqcup}\hspace{-0.07cm}\boldsymbol{\sqcup}}}}
\newcommand{\shuff}{\ensuremath{\scriptscriptstyle{\boldsymbol{\sqcup}\hspace{-0.07cm}\boldsymbol{\sqcup}}}\displaystyle{}}

\def\A{{\mathcal {A}}}
\def\AW{{\mathcal {A}^{\ast}}}

\def\NN{{\mathbb {N}}}
\def\QQ{{\mathbb {Q}}}
\def\RR{{\mathbb {R}}}

\DeclareMathOperator{\Li}{Li}

\advance\textheight by \topskip

\newtheorem{coroll}{Corollary}
\newtheorem{theorem}{Theorem}
\newtheorem{lemma}{Lemma}

\theoremstyle{definition}

\title{On a reciprocity law for finite multiple zeta values}

\author[M.~Kuba]{Markus Kuba}

\address{Markus Kuba\\
Institut f{\"u}r Diskrete Mathematik und Geometrie\\
Technische Universit\"at Wien\\
Wiedner Hauptstr. 8-10/104\\
1040 Wien, Austria} %
\email{kuba@dmg.tuwien.ac.at}

\author[H.~Prodinger]{Helmut Prodinger}
\address{Helmut Prodinger\\
Department of Mathematics\\
University of Stellenbosch\\
7602 Stellenbosch\\
South Africa}
\email{hproding@sun.ac.za}

\begin{document}

\begin{abstract}
It was shown in~\cite{KirProd98,ProSchnKu} that harmonic numbers satisfy certain reciprocity relations, which are in particular useful for the analysis of the quickselect algorithm. The aim of this work is to show that a reciprocity relation from~\cite{KirProd98,ProSchnKu} can be generalized to finite variants of multiple zeta values, involving a finite variant of the shuffle identity for multiple zeta values. We present the generalized reciprocity relation and furthermore a simple elementary proof of the shuffle identity using only partial fraction decomposition.
We also present an extension of the reciprocity relation to weighted sums.
\end{abstract}

\thanks{The author M.~K.~was supported by the Austrian Science
Foundation FWF, grant S9608-N13}

\maketitle

\section{Introduction}
Let $H_n=\sum_{k=1}^{n}1/k$ denote the $n$-th harmonic number and $H_{n}^{(s)}=\sum_{k=1}^{n}1/k^{s}$ the $n$-th harmonic number of order $s$, with $n,s\in\NN$ and $H_n=H_n^{(1)}$. Kirschenhofer and Prodinger~\cite{KirProd98} analysed the variance
of the number of comparisons of the quickselect algorithm~\cite{HoFIND}, and derived a reciprocity relation for (first order) harmonic
numbers. Subsequently, the reciprocity relation of~\cite{KirProd98} was generalized in~\cite{ProSchnKu}, where the following identity was derived.
\begin{equation}
\label{Pro1}
\begin{split}
    \sum_{k=1}^{j}\frac{H_{N-k}^{(a)}}{k^b}%
    + \sum_{k=1}^{N+1-j}\frac{H_{N-k}^{(b)}}{k^a}
    &   = -\frac{1}{j^b(N+1-j)^a}+ H^{(b)}_{j}H^{(a)}_{N+1-j} + R_N^{(a,b)},
\end{split}
\end{equation}
where $R^{(a,b)}_N=\sum_{k=1}^{N}\frac{H_{N-k}^{(a)}}{k^b}$, which can be evaluated into a finite analog of the so-called
Euler identity for $\zeta(a)\zeta(b)$ stated below,
\begin{equation}
\label{Rn}
\begin{split}
    R^{(a,b)}_N & =   \sum_{i=1}^a\binom{i+b-2}{b-1}\zeta_N(i+b-1,a+1-i)
        + \sum_{i=1}^b\binom{i+a-2}{a-1} \zeta_N(i+a-1,b+1-i),
\end{split}
\end{equation}
where the multiple zeta values~\cite{Bor0,Bor,Bow,Bow2,Bow3,Zag}, and its finite counterpart are defined as follows:
\begin{align*}
     \zeta(\mathbf{a}) =\zeta(a_1,\dots,a_r) &:= \sum_{ n_1 > n_2 > \dots > n_{r}\geq
    1} \frac{1}{ n_1^{a_1} n_2^{a_2} \dots n_r^{a_r} },\\
     \zeta_N(\mathbf{a}) =\zeta_N(a_1,\dots,a_r) &:= \sum_{ N \geq n_1 > n_2 > \dots > n_{r}\geq
    1} \frac{1}{ n_1^{a_1} n_2^{a_2} \dots n_r^{a_r} }.
\end{align*}
Note that $\zeta_N(a)=H_N^{(a)}$. Let $w=\sum_{i=1}^{r}a_i$ denote
the weight and $d=r$ the depth of (finite) multiple
zeta values. The aim of this note is to derive a generalization of
the reciprocity relation~\eqref{Pro1}, stated below in
Theorem~\ref{Prothe1}, by considering the more general sums
\begin{equation*}
\displaystyle{\sum_{k=1}^{j}\frac{\zeta_{k-1}(b_2,\dots,b_s)\zeta_{N-k}(a_1,\dots,a_r)}{k^{b_1}}}+\displaystyle{\sum_{k=1}^{N+1-j}\frac{\zeta_{k-1}(a_2,\dots,a_r)\zeta_{N-k}(b_1,\dots,b_r)}{k^{a_1}}},
\end{equation*}
instead of the previously considered sums $\sum_{k=1}^{j}\frac{H_{N-k}^{(a)}}{k^b}$ and $\sum_{k=1}^{N+1-j}\frac{H_{N-k}^{(b)}}{k^a}$.
Our generalization involves a finite variant of the shuffle identity
for multiple zeta values, for which we give an elementary proof
using partial fraction decomposition. We discuss the close relation between this finite variant of the shuffle identity and the shuffle identity for generalized polylogarithm functions; we will see that the finite variant of the shuffle identity is equivalent to the shuffle identity for generalized polylogarithm functions. To simplify the presentation of this work we will frequently use the shorthand notations $\mathbf{a}=(a_1,\dots,a_r)$, $\mathbf{a}_2=(a_2,\dots,a_r)$ and $\mathbf{b}=(b_1,\dots,b_s)$, $\mathbf{b}_2=(b_2,\dots,b_s)$, respectively, with $r,s\in\NN$ and $a_i,b_k\in\NN$ for $1\le i\le r$ and $1\le k\le s$.
%Consequently,
%$\sum_{k=1}^{j}\frac{\zeta_{k-1}(b_2,\dots,b_s)\zeta_{N-k}(a_1,\dots,a_r)}{k^{b_1}}=\sum_{k=1}^{j}\frac{\zeta_{k-1}(\mathbf{b}_2)\zeta_{N-k}(\mathbf{a})}%{k^{b_1}}$, and %$\sum_{k=1}^{N+1-j}\frac{\zeta_{k-1}(a_2,\dots,a_r)\zeta_{N-k}(b_1,\dots,b_r)}{k^{a_1}}=\sum_{k=1}^{N+1-j}\frac{\zeta_{k-1}(\mathbf{a}_2)\zeta_{N-k}(\mathbf{b})%}{k^{a_1}}$.

\section{The reciprocity relation for finite multiple zeta values}
We will state our main theorem below, and subsequently discuss its proof and the precise definition of the shuffle relation for multiple zeta values.
\begin{theorem}
The multiple zeta values $\zeta_N(\mathbf{a})=\zeta_N(a_1,\dots,a_r)$ and $\zeta_N(\mathbf{b})=\zeta_N(b_1,\dots,b_s)$ satisfy the following reciprocity relation.
\label{Prothe1}
\begin{equation*}
\begin{split}
&\sum_{k=1}^{j}\frac{\zeta_{k-1}(b_2,\dots,b_s)\zeta_{N-k}(a_1,\dots,a_r)}{k^{b_1}}
+ \sum_{k=1}^{N+1-j}\frac{\zeta_{k-1}(a_2,\dots,a_r)\zeta_{N-k}(b_1,\dots,b_r)}{k^{a_1}}\\
&\quad=\zeta_{N+1-j}(\mathbf{a})\zeta_{j}(\mathbf{b})-\frac{\zeta_{j-1}(\mathbf{b}_2)\zeta_{N-j}(\mathbf{a}_2)}{j^{b_1}(N+1-j)^{a_1}}+
R_N(\mathbf{a};\mathbf{b}).
\end{split}
\end{equation*}
Here $R_N(\mathbf{a};\mathbf{b})=\sum_{k=1}^{N}\frac{\zeta_{N-k}(\mathbf{b})\zeta_{k-1}(a_2,\dots,a_r)}{k^{a_1}}=R_N(\mathbf{b};\mathbf{a})$ satisfies a finite counterpart of the shuffle identity $\zeta(\mathbf{a})\zeta(\mathbf{b})=\zeta(\mathbf{a}\shuff\mathbf{b})$ for the multiple zeta value, $R_N(\mathbf{a};\mathbf{b})=\zeta_N(\mathbf{a}\shuff\mathbf{b})$.
\end{theorem}

\begin{coroll}
We obtain the complementary identity
\begin{equation*}
\begin{split}
&\sum_{k=1}^{j-1}\frac{\zeta_{k}(\mathbf{b})\zeta_{N-k-1}(\mathbf{a}_2)}{(N-k)^{a_1}}+
\sum_{k=1}^{N-j}\frac{\zeta_{k}(\mathbf{a})\zeta_{N-k-1}(\mathbf{b}_2)}{(N-k)^{b_1}}
=\frac{\zeta_{j-1}(\mathbf{b})\zeta_{N-j}(\mathbf{a}_2)}{(N+1-j)^{a_1}}+\frac{\zeta_{N-j}(\mathbf{a})\zeta_{j-1}(\mathbf{b}_2)}{j^{b_1}}\\
&\quad-\zeta_{N+1-j}(\mathbf{a})\zeta_{j}(\mathbf{b})+\frac{\zeta_{j-1}(\mathbf{b}_2)\zeta_{n-j}(\mathbf{a}_2)}{j^{b_1}(N+1-j)^{a_1}}+
R_N(\mathbf{a};\mathbf{b}).
\end{split}
\end{equation*}
\end{coroll}

Next we state an immediate asymptotic implication of our previous result.
\begin{coroll}
\label{SHUFFcor1}
For $N=2n+1$, $j=n+1$, with  $a_1,b_1\in\NN\setminus\{1\}$ and $n\to\infty$ we obtain the following result.
\begin{equation*}
\begin{split}
&\lim_{n\to\infty}\biggl(\sum_{k=1}^{j}\frac{\zeta_{k-1}(\mathbf{b}_2)\zeta_{N-k}(\mathbf{a})}{k^{b_1}}+\sum_{k=1}^{N+1-j}\frac{\zeta_{k-1}(\mathbf{a}_2)\zeta_{N-k}(\mathbf{b})}{k^{a_1}}\biggr)=2\zeta(\mathbf{a})\zeta(\mathbf{b}).
\end{split}
\end{equation*}
\end{coroll}
In order to prove Theorem~\ref{Prothe1} we proceed as follows.
\begin{equation*}
\begin{split}
\sum_{k=1}^{j}\frac{\zeta_{k-1}(\mathbf{b}_2)\zeta_{N-k}(\mathbf{a})}{k^{b_1}} & = \zeta_{N-j}(\mathbf{a})\zeta_{j}(\mathbf{b})
     +\sum_{k=1}^{j}\frac{\zeta_{k-1}(\mathbf{b}_2)}{k^{b_1}}\sum_{\ell=N+1-j}^{N-k} \frac{\zeta_{\ell-1}(\mathbf{a}_2)}{\ell^{a_1}} \\
    &= \zeta_{N-j}(\mathbf{a})\zeta_{j}(\mathbf{b}) + \sum_{l=N+1-j}^{N-1}\frac{\zeta_{\ell-1}(\mathbf{a}_2)}{\ell^{a_1}} \sum_{k=1}^{N-l}\frac{\zeta_{k-1}(\mathbf{b}_2)}{k^{b_1}}\\
    & = \zeta_{N-j}(\mathbf{a})\zeta_{j}(\mathbf{b})+ \sum_{\ell=N+1-j}^{N-1}\frac{\zeta_{\ell-1}(\mathbf{a}_2)\zeta_{N-\ell}(\mathbf{b})}{\ell^{a_1}}\\
   & = \zeta_{N+1-j}(\mathbf{a})\zeta_{j}(\mathbf{b})+  \frac{\zeta_{N-j}(\mathbf{a}_2)\zeta_{j-1}(\mathbf{b}_2)}{(N+1-j)^{a_1}j^{b_1}} +\sum_{\ell=N+2-j}^{N}\frac{\zeta_{\ell-1}(\mathbf{a}_2)\zeta_{N-\ell}(\mathbf{b})}{\ell^{a_1}}.
\end{split}
\end{equation*}
This proves the first part of Theorem~\ref{Prothe1} and
\begin{equation*}
R_N(\mathbf{a};\mathbf{b})=\sum_{k=1}^{N}\frac{\zeta_{N-k}(\mathbf{b})\zeta_{k-1}(a_2,\dots,a_r)}{k^{a_1}}.
\end{equation*}
For the evaluation of $R_N(\mathbf{a};\mathbf{b})$ we note that $R_{0}(\mathbf{a};\mathbf{b})=0$, and further
\begin{equation}
\label{SHUFFrn}
R_N(\mathbf{a};\mathbf{b})=\sum_{k=1}^{N}\big(R_{k}(\mathbf{a};\mathbf{b})-R_{k-1}(\mathbf{a};\mathbf{b})\big).
\end{equation}
We have
\begin{equation*}
R_N(\mathbf{a};\mathbf{b})-R_{N-1}(\mathbf{a};\mathbf{b})=\sum_{k=1}^{N-1}\frac{\zeta_{k-1}(a_2,\dots,a_r)\zeta_{N-1-k}(b_2,\dots,b_s)}{(N-k)^{b_1}k^{a_1}}.
\end{equation*}
Now we use partial fraction decomposition\footnote{This identity has been rediscovered many times. For a fascinating historic account, see~\cite{KoornSchloss2008}.}, which appears already in~\cite{Niel65},
\begin{align}
\label{frac}
\frac1{k^a(N-k)^b}=\sum_{i=1}^a\frac{\binom{i+b-2}{b-1}}{N^{i+b-1}k^{a+1-i}}+
\sum_{i=1}^b\frac{\binom{i+a-2}{a-1}}{N^{i+a-1}(N-k)^{b+1-i}},
\end{align}
and obtain
\begin{equation*}
\begin{split}
\sum_{k=1}^{N-1}\frac{\zeta_{k-1}(a_2,\dots,a_r)\zeta_{N-1-k}(b_2,\dots,b_s)}{(N-k)^{b_1}k^{a_1}}&=
\sum_{i=1}^{a_1}\sum_{k=1}^{N-1}\frac{\binom{i+b_1-2}{b_1-1}\zeta_{k-1}(a_2,\dots,a_r)\zeta_{N-1-k}(b_2,\dots,b_s)}{N^{i+b_1-1}k^{a_1+1-i}}\\
&+ \sum_{i=1}^{b_1}\sum_{k=1}^{N-1}\frac{\binom{i+a_1-2}{a_1-1}\zeta_{k-1}(a_2,\dots,a_r)\zeta_{N-1-k}(b_2,\dots,b_s)}{N^{i+a_1-1}(N-k)^{b_1+1-i}}.
\end{split}
\end{equation*}
Consequently, by summing up according to~\eqref{SHUFFrn} we get the following recurrence relation for $R_N(\mathbf{a};\mathbf{b})$.
\begin{equation}
\begin{split}
\label{SHUFFeqnT1}
R_N(\mathbf{a};\mathbf{b})&=\sum_{i=1}^{a_1}\sum_{n_1=1}^{N}\frac{\binom{i+b_1-2}{b_1-1}}{n_1^{i+b_1-1}}
R_{n_1-1}(a_1+1-i,a_2,\dots,a_r;b_2,\dots,b_s) \\
&+\sum_{i=1}^{b_1}\sum_{n_1=1}^{N}\frac{\binom{i+a_1-2}{a_1-1}}{n_1^{i+a_1-1}}R_{n_1-1}(a_2,\dots,a_r;b_1+1-i,b_2,\dots,b_s).
\end{split}
\end{equation}
This recurrence relation suggests that there exists an evaluation of $R_N(\mathbf{a};\mathbf{b})$ into sums of finite
multiple zeta values, all of them having weight $w=\sum_{i=1}^{r}a_r+\sum_{i=1}^{s}b_i$ and depth $d=r+s$.
In order to specify this evaluation we need to introduce the shuffle algebra for (finite) multiple zeta values.

\subsection{The shuffle algebra}
Let $\A$ denote a finite non-commutative alphabet consisting of a set of letters. A word $\mathbf{w}$ on the alphabet $\A$
consists of a sequence of letters from $\A$. Let $\AW$ denote the
set of all words on the alphabet $\A$. A polynomial on $\A$ over $\QQ$ is a rational linear
combination of words on $\A$. The set of all such polynomials is denoted by $\QQ\langle\A\rangle$. Let the \emph{shuffle product} $\shuff$ be defined on $\QQ\langle\A\rangle$ as follows: for any $\mathbf{w},\mathbf{v}\in\AW$ with $\mathbf{w}=x_1\dots x_n$, $\mathbf{v}=x_{n+1}\dots x_{n+m}$, $x_i\in\A$ for $1\le i \le n+m$
\begin{equation}
\mathbf{w} \shuff \mathbf{v} := \sum x_{\sigma(1)} x_{\sigma(2)}\dots  x_{\sigma(n+m)},
\end{equation}
where the sum runs over all $\binom{n+m}{n}$ permutations $\sigma\in\mathfrak{S}_{n+m}$ which satisfy $\sigma^{-1}(j)<\sigma^{-1}(k)$ for all $1\le j< k\le n$ and $n+1\le j <k \le n+m$. The sum is over all words of length $n+m$, counting multiplicities,
in which the relative orders of the letters $x_1,\dots,x_n$ and $x_{n+1},\dots,x_{n+m}$
are preserved. The term ``shuffle'' is used because such permutations arise in riffle shuffling a deck of $n + m$
cards cut into one pile of $n$ cards and a second pile of $m$ cards~\cite{Bow3}.
Equivalently, we can  recursively define the shuffle product as follows.
\begin{align}
\label{SHUFFrec}
&\forall \mathbf{w}\in\AW, &\epsilon \shuff \mathbf{w} &= \mathbf{w} \shuff \epsilon =\mathbf{w},\nonumber\\
&\forall x,y\in \A,\quad \mathbf{w},\mathbf{v}\in\AW, & x\mathbf{w}\shuff y\mathbf{v}&=
x(\mathbf{w}\shuff y\mathbf{v})+y(x\mathbf{w}\shuff \mathbf{v}).
\end{align}

\subsection{The shuffle algebra and multiple zeta values}
Let $\mathbf{a}$ and $\mathbf{b}$ denote the multi-indices
$\mathbf{a}=(a_1,\dots,a_r)$ and $\mathbf{b}=(b_1,\dots,b_s)$ with
$a_i,b_j\in\NN$ for $1\le i \le r$, $1\le j \le s$. To any
multi-index we associate a unique word over the non commutative
alphabet $\mathcal{A}=\{\omega_0,\omega_1\}$. Let $A=A(\mathbf{a})$,
$B=B(\mathbf{b})$ such that
$A:=\omega_0^{a_1-1}\omega_1\omega_0^{a_2-1}\omega_1\dots\omega_0^{a_r-1}\omega_1$
and
$B:=\omega_0^{b_1-1}\omega_1\omega_0^{b_2-1}\omega_1\dots\omega_0^{b_s-1}\omega_1$.
To each word we associate a finite multiple zeta values by the
following linear correspondence:
$Z_N(\omega_1\omega_0^{a-1})=\zeta_N(a)$, and in general
\begin{equation}
\begin{split}
\label{SHUFFcorres}
Z_N\Big(\omega_0^{a_1-1}\omega_1\omega_0^{a_2-1}\omega_1\dots\omega_0^{a_r-1}\omega_1\Big)
=\sum_{n_1=1}^{N}\frac{1}{n_1^{a_1}}Z_{n_1-1}(\omega_0^{a_2-1}\omega_1\dots\omega_0^{a_r-1}\omega_1)=\zeta_N(a_1,\dots,a_r).
\end{split}
\end{equation}
Moreover, assuming that $\mathbf{a}_{\ell}=(a_{\ell,1},\dots,a_{\ell,r_{\ell}})$ and with $1\le \ell \le h$, $q_\ell\in\QQ$, and $A_{\ell}=A(\mathbf{a}_{\ell})$, we get by linearity of the correspondence \begin{equation*}
Z_N\biggl(\sum_{\ell=1}^{h}q_{\ell}A_{\ell}\biggr)=\sum_{\ell=1}^{h}q_{\ell}Z_N(A_{\ell})=\sum_{\ell=1}^{h}q_{\ell}\zeta_N(a_{\ell,1},\dots,a_{\ell,r_{\ell}}).
\end{equation*}

\smallskip

We observe that the partial fraction decomposition~\eqref{frac} of $\frac1{k^a(N-k)^b}$ above mimics the basic shuffle identity for words $A=\omega_0^{a-1}\omega_1$, $B=\omega_0^{b-1}\omega_1$,
\begin{equation*}
A\shuff B=\sum_{i=0}^{a-1}\binom{b-1+i}{b-1}\omega_0^{b-1+i}\omega_1\omega_0^{a-1-i}\omega_1 +\sum_{i=0}^{b_1-1}\binom{a-1+i}{a-1}\omega_0^{a-1+i}\omega_1\omega_0^{b-1-i}\omega_1,
\end{equation*}
which appeared in Hoang and Petitot~\cite{Minh}.

\smallskip

The key to the explicit evaluation of $R_N(\mathbf{a};\mathbf{b})$ is the following result
concerning the shuffling of the words $A$ and $B$, associated to the multi-indices $\mathbf{a}$ and $\mathbf{b}$.
\begin{lemma}
\label{SHUFFlem1}
Let $A:=\omega_0^{a_1-1}\omega_1\omega_0^{a_2-1}\omega_1\dots\omega_0^{a_r-1}\omega_1$ and
$B:=\omega_0^{b_1-1}\omega_1\omega_0^{b_2-1}\omega_1\dots\omega_0^{b_s-1}\omega_1$,
with $\mathbf{a}=(a_1,\dots,a_r)$ and $(b_1,\dots,b_s)$ with $a_i,b_j\in\NN$, $1\le i \le r$, $1\le j \le s$.
We have
\begin{equation*}
\begin{split}
A\shuff B&= %\sum_{i=1}^{a_1}\binom{i+b_1-2}{b_1-1}\omega_0^{i+b_1-2}\omega_1 \big(
%\omega_0^{a_1-i}\omega_1\omega_0^{a_2-1}\omega_1\dots\omega_0^{a_r-1}\omega_1 \shuff \omega_0^{b_2-1}\omega_1\dots\omega_0^{b_s-1}\omega_1\big)\\
%&+\sum_{i=1}^{b_1}\binom{i+a_1-2}{a_1-1}\omega_0^{i+a_1-2}\omega_1 \big(
%\omega_0^{a_2-1}\omega_1\dots\omega_0^{a_r-1}\omega_1 \shuff \omega_0^{b_1-1-i}\omega_1\omega_0^{b_2-1}\omega_1\dots\omega_0^{b_s-1}\omega_1\big)\\
%&=
\sum_{i=1}^{a_1}\binom{i+b_1-2}{b_1-1}\omega_0^{i+b_1-2}\omega_1(A_i'\shuff B_2)
+ \sum_{i=1}^{b_1}\binom{i+a_1-2}{a_1-1}\omega_0^{i+a_1-2}\omega_1(A_2\shuff B_i'),
\end{split}
\end{equation*}
with $A_i':=\omega_0^{a_1-i}\omega_1\omega_0^{a_2-1}\omega_1\dots\omega_0^{a_r-1}\omega_1$,
$B_i':=\omega_0^{b_1-i}\omega_1\omega_0^{b_2-1}\omega_1\dots\omega_0^{b_s-1}\omega_1$ and further
$A_2:=\omega_0^{a_2}\omega_1\dots\omega_0^{a_r-1}\omega_1$,
$B_2:=\omega_0^{b_2-1}\omega_1\dots\omega_0^{b_s-1}\omega_1$.
\end{lemma}
The special case $r=s=1$ is a result of Hoang and Petitot~\cite{Minh}; we simply use the recursive definition
of the shuffle product~\eqref{SHUFFrec} and obtain the result of Lemma~\ref{SHUFFlem1}.

Now we are ready to prove the evaluation of $R_N(\mathbf{a};\mathbf{b})$.
Let $A=A(\mathbf{a})$ and $B=A(\mathbf{b})$ denote the words associated to the multi-indices $\mathbf{a}$ and $\mathbf{b}$,
\begin{equation}
\label{SHUFFfinite}
R_N(\mathbf{a};\mathbf{b})=Z_N(A\shuff B)=: \zeta_N(\mathbf{a}\shuff \mathbf{b}).
\end{equation}
We use induction with respect to the depth $d=r+s$. The result clearly holds for $d=2$, \eqref{Rn}, as shown in~\cite{ProSchnKu}.
Assuming the result for all depth $r+s<d$ we obtain according to the induction hypothesis
\begin{equation*}
\begin{split}
R_N(\mathbf{a};\mathbf{b})&=\sum_{i=1}^{a_1}\sum_{n_1=1}^{N}\frac{\binom{i+b_1-2}{b_1-1}}{n_1^{i+b_1-1}}
R_{n_1-1}(a_1+1-i,a_2,\dots,a_r;b_2,\dots,b_s) \\
&\quad+\sum_{i=1}^{b_1}\sum_{n_1=1}^{N}\frac{\binom{i+a_1-2}{a_1-1}}{n_1^{i+a_1-1}}R_{n_1-1}(a_2,\dots,a_r;b_1+1-i,b_2,\dots,b_s)\\
&=\sum_{i=1}^{a_1}\sum_{n_1=1}^{N}\frac{\binom{i+b_1-2}{b_1-1}}{n_1^{i+b_1-1}}
Z_{n_1-1}\big(A_i\shuff B_2 \big)+\sum_{i=1}^{b_1}\sum_{n_1=1}^{N}\frac{\binom{i+a_1-2}{a_1-1}}{n_1^{i+a_1-1}}Z_{n_1-1}\big( A_2\shuff B_i\big).
\end{split}
\end{equation*}
By~\eqref{SHUFFcorres} and Lemma~\ref{SHUFFlem1}, using the notations of Lemma~\ref{SHUFFlem1}, we get 
\begin{equation*}
\begin{split}
Z_N\big(A\shuff B\big)=\sum_{i=1}^{a_1}\sum_{n_1=1}^{N}\frac{\binom{i+b_1-2}{b_1-1}}{n_1^{i+b_1-1}}
Z_{n_1-1}\big(A_i\shuff B_2 \big)+\sum_{i=1}^{b_1}\sum_{n_1=1}^{N}\frac{\binom{i+a_1-2}{a_1-1}}{n_1^{i+a_1-1}}Z_{n_1-1}\big( A_2\shuff B_i\big).
\end{split}
\end{equation*}
Consequently, 
\begin{equation}
R_N(\mathbf{a};\mathbf{b})=Z_N(A\shuff B)=\zeta_N(\mathbf{a}\shuff \mathbf{b}).
\end{equation}
This proves the stated result for $R_N(\mathbf{a};\mathbf{b})$ and the remaining part of Theorem~\ref{Prothe1}.
Moreover we observe that the depths $d=r+s$ and the weights $w=\sum_{i=1}^{r}a_i+\sum_{k=1}^{s}b_k$ of the finite multiple zeta values are all the same. Corollary~1 can easily be deduced by noting that the sum of the left hand sides of Corollary~1 and Theorem~\ref{Prothe1}
add up to $R_n(\mathbf{a};\mathbf{b})$ with respect to two extra terms. Now we turn to the proof of Corollary~\ref{SHUFFcor1}. Since for $N=2n+1$ and $j=n+1$ and $n\to\infty$ we have
\begin{equation*}
\begin{split}
&\lim_{n\to\infty}\zeta_{j}(\mathbf{b})\zeta_{N+1-j}(\mathbf{b})=\lim_{n\to\infty} \zeta_{n+1}(\mathbf{b})\zeta_{n+1}(\mathbf{b})=\zeta(\mathbf{a})\zeta(\mathbf{b}),\\
&\lim_{n\to\infty}\frac{\zeta_{n}(b_2,\dots,b_s)\zeta_{n}(a_2,\dots,a_r)}{(n+1)^{a_1+b_1}}=0,\\
&\lim_{n\to\infty}R_N(\mathbf{a};\mathbf{b})=\lim_{n\to\infty}\zeta_{2n+1}\big(\mathbf{a}\shuff \mathbf{b}\big)=\zeta\big(\mathbf{a}\shuff \mathbf{b}\big)=\zeta(\mathbf{a})\zeta(\mathbf{b}),
\end{split}
\end{equation*}
the result of Corollary~\ref{SHUFFcor1} can be immediately deduced from Theorem~\ref{Prothe1}. Note that the last identity is the well known
shuffle identity for multiple zeta values; we refer the reader to the excellent article~\cite{Bor}.

\section{Polylogarithms and the finite shuffle identity}
Let $\Li_{\mathbf{a}}(z)=\Li_{a_1,\dots,a_r}(z)$ denote the (multiple) polylogarithm function with parameters $a_1,\dots,a_r$,
defined by
\begin{equation}
\label{SHUFFpolylog}
\Li_{\mathbf{a}}(z)=\Li_{a_1,\dots,a_r}(z)=\sum_{n_1>n_2>\dots>n_r\ge 1}\frac{z^{n_1}}{n_1^{a_1}n_2^{a_2}\dots n_{r}^{a_r}},
\end{equation}
We note that the value $R_N(\mathbf{a};\mathbf{b})$ can be obtained in the following way.
\begin{equation*}
R_N(\mathbf{a};\mathbf{b})=\sum_{k=1}^{N}\frac{\zeta_{k-1}(a_2,\dots,a_r)\zeta_{N-k}(b_1,\dots,b_s)}{k^{a_1}}=[z^N]\frac{\Li_{\mathbf{a}}(z)\Li_{\mathbf{b}}(z)}{1-z}
\end{equation*}
Consequently, our finite shuffle identity~\eqref{SHUFFfinite} for $R_N(\mathbf{a};\mathbf{b})$ is equivalent to the following shuffle identity for polylogarithm functions.
\begin{equation*}
\Li_{\mathbf{a}}(z)\Li_{\mathbf{b}}(z)=\Li_{\mathbf{a}\shuffk\mathbf{b}}(z).
\end{equation*}
Note that by evaluating at $z=1$ the shuffle identity for polylogarithm functions implies the shuffle identity for multiple zeta values.
The identity above is well known, see for example the article~\cite{Bor}. The shuffle identity for polylogarithm functions is due to the iterated Drinfeld integral representation of polylogarithm functions and multiple zeta values due to Kontsevich~\cite{Zag}. As remarked in~\cite{Bor} the shuffle product holds since the product of two simplex integrals consists of a sum of simplex integrals over all possible interlacings of the respective variables of integration. Our result for $R_N(\mathbf{a};\mathbf{b})$ implies that the shuffle for polylogarithm
functions, and also for multiple zeta values, can be derived using only basic partial fraction decomposition.

\section{The reciprocity relation for weighted multiple zeta values}
Results similar to Theorem~\ref{Prothe1} and Corollary~\ref{SHUFFcor1} can be obtained for products of weighted finite multiple zeta values,
$\zeta_N(a_1,a_2,\dots,a_r;\sigma_1,\dots,\sigma_r)$, $\sigma_i\in\RR\setminus\{0\}$ for $1\le i\le r$, defined as follows:
\begin{align*}
   \zeta_{N}(\mathbf{a},\boldsymbol{\sigma})= \zeta_{N}(a_1,a_2,\dots,a_r;\sigma_1,\dots,\sigma_r) &= \sum_{ N\ge n_1 > n_2 > \dots >
    n_{r} \ge 1}\frac{1}{\prod_{i=1}^{r}n_i^{a_i}\sigma_i^{n_i}}.
\end{align*}
Of particular interest are the cases $\sigma_i\in\{\pm 1\}$ corresponding to a mixture of alternating and non-alternating signs.
We only state the result generalizing Theorem~\ref{Prothe1}, with respect to the notations $\mathbf{a}_2=(a_2,\dots,a_r)$, $\boldsymbol{\sigma}_2=(\sigma_2,\dots,\sigma_r)$ and the corresponding notations for $\mathbf{b}_2$ and $\boldsymbol{\tau}_2$,
and leave the generalizations of Corollaries 1,2 to the reader.
\begin{theorem}
The multiple zeta values $\zeta_N(\mathbf{a},\boldsymbol{\sigma})$ and $\zeta_N(\mathbf{b},\boldsymbol{\tau})$ with weights $\boldsymbol{\sigma}$ and $\boldsymbol{\tau}$ satisfy the following reciprocity relation.
\begin{equation*}
\begin{split}
&\sum_{k=1}^{j}\frac{\zeta_{k-1}(\mathbf{b}_2,\boldsymbol{\tau}_2)\zeta_{N-k}(\mathbf{a},\boldsymbol{\sigma})}{k^{b_1}\tau_1^{k}}
+ \sum_{k=1}^{N+1-j}\frac{\zeta_{k-1}(\mathbf{a}_2,\boldsymbol{\sigma}_2)\zeta_{N-k}(\mathbf{b},\boldsymbol{\tau})}{k^{a_1}\sigma_1^{k}}\\
&\quad=\zeta_{N+1-j}(\mathbf{a},\boldsymbol{\sigma})\zeta_{j}(\mathbf{b},\boldsymbol{\tau})-\frac{\zeta_{j-1}(\mathbf{b}_2,\boldsymbol{\tau}_2)\zeta_{N-j}(\mathbf{a}_2,\boldsymbol{\sigma}_2)}{\tau_1^{j}j^{b_1}\sigma_1^{N+1-j}(N+1-j)^{a_1}}+
R_N(\mathbf{a},\boldsymbol{\sigma};\mathbf{b},\boldsymbol{\tau}).
\end{split}
\end{equation*}
Here $R_N(\mathbf{a},\boldsymbol{\sigma};\mathbf{b},\boldsymbol{\tau})=\sum_{k=1}^{N}\frac{\zeta_{N-k}(\mathbf{b},\boldsymbol{\tau})\zeta_{k-1}(\mathbf{a}_2,\boldsymbol{\sigma}_2)}{\sigma_1^{k}k^{a_1}}=R_N(\mathbf{b},\boldsymbol{\tau};\mathbf{a},\boldsymbol{\sigma})$ satisfies an analogue of the shuffle identity with respect to the weights $\boldsymbol{\sigma}$ and $\boldsymbol{\tau}$.
\end{theorem}

The proof of Theorem~\ref{Prothe1} can easily be adapted to the weighted case. Hence, we only elaborate on the main new difficulty, namely the evaluation of the quantity
\begin{equation*}
R_N(\mathbf{a},\boldsymbol{\sigma};\mathbf{b},\boldsymbol{\tau})=\sum_{k=1}^{N}\frac{\zeta_{N-k}(\mathbf{b},\boldsymbol{\tau})\zeta_{k-1}(\mathbf{a}_2;\boldsymbol{\sigma}_2)}{\sigma_1^{k}k^{a_1}}.
\end{equation*}
Proceeding as before, i.e.~taking differences and using partial fraction decomposition,
we obtain the recurrence relation
\begin{equation*}
\begin{split}
R_N(\mathbf{a},\boldsymbol{\sigma};\mathbf{b},\boldsymbol{\tau})&=\sum_{i=1}^{a_1}\sum_{n_1=1}^{N}\frac{\binom{i+b_1-2}{b_1-1}}{n_1^{i+b_1-1}\tau_1^{n_1}}R_{n_1-1}(a_1+1-i,\mathbf{a}_2,\frac{\tau_1}{\sigma_1},\boldsymbol{\sigma}_2 ; \mathbf{b}_2,\boldsymbol{\tau}_2) \\
&+\sum_{i=1}^{b_1}\sum_{n_1=1}^{N}\frac{\binom{i+a_1-2}{a_1-1}}{n_1^{i+a_1-1}\sigma_1^{n_1}}R_{n_1-1}(\mathbf{a}_2,\boldsymbol{\sigma}_2;b_1+1-i,\mathbf{b}_2,\frac{\sigma_1}{\tau_1},\boldsymbol{\tau}_2).
\end{split}
\end{equation*}
Consequently, the value $R_N(\mathbf{a},\boldsymbol{\sigma};\mathbf{b},\boldsymbol{\tau})$ can be evaluated into
sums of weighted finite multiple zeta values according to a shuffle identity with respect to the weights $\boldsymbol{\sigma}$ and $\boldsymbol{\tau}$.
We omit the precise definition of this generalization and leave the details to the interested reader.

\section*{Conclusion}
We presented a reciprocity relation for finite multiple zeta values, extending the previous results of~\cite{KirProd98,ProSchnKu}.
The reciprocity relation involves a shuffle product identity for (finite) multiple zeta values, for which we gave a simple proof using only partial fraction decomposition. Moreover, we also presented the reciprocity relation for weighted finite multiple zeta values.

\end{document}